\documentclass{elsarticle}
\input{packages.sty}
\input intmacros.sty

\usepackage{todonotes}

\newcommand{\ignore}[1]{}

\begin{document}

\title{The Laurent-Horner method for validated evaluation of Chebyshev expansions}

\author[ja]{Jared L. Aurentz}
\ead{jared.aurentz@icmat.es}

\author[bh]{Behnam Hashemi\corref{cor1}}
\ead{hoseynhashemi@gmail.com}

\cortext[cor1]{Corresponding author}

\address[ja]{Instituto de Ciencias Matem\'aticas, Universidad Aut\'onoma de Madrid, Madrid, Spain}


\address[bh]{Department of Mathematics, College of Sciences, Shiraz University, Shiraz 71454, Iran}

\begin{abstract}
We develop a simple two-step algorithm for enclosing Chebyshev expansions whose cost is linear in terms of the polynomial degree.\ The algorithm first transforms the expansion from Chebyshev to the Laurent basis and then applies the interval Horner method. It outperforms the existing eigenvalue-based methods if the degree is high or the evaluation point is close to the boundaries of the domain. 
\end{abstract}


\maketitle

\section{Introduction}
\label{intro}Let 
\begin{equation}
\label{ChebExpand: eq}
p(x) = \sum_{k=0}^{n}c_k T_k(x)
\end{equation}
be a finite Chebyshev expansion in which 
\[T_k(x) = \cos(k \cos^{-1}(x))\]
is the $k$-th Chebyshev polynomial of the first kind defined on $[-1, 1]$, and $c_k$ is the $k$-th constant real coefficient.\ Recently, several techniques are developed in \cite{Hashemi19}\ignore{and compared with the older ones in order 
to compute narrow enclosures for the value of the polynomial $p$.} for validated evaluation of $p$. We refer the reader to \cite{Alefeld83,Neumaier90,Mayer17} for the basics of interval arithmetic, in particular the enclosure property, dependency issue, the wrapping effect, and how directed roundings are employed with floating point arithmetic to efficiently compute results which are guaranteed to be mathematically correct. Enclosing Chebyshev expansions has application in computer-assisted existence proofs of the spherical t-designs \cite{Chen11}, ultra arithmetic\ignore{removed \cite{Epstein82} to free space!}, Chebyshev models \cite{Joldes11} and automatic a posteriori forward error analysis of floating point evaluation of Chebyshev expansions; see \cite{Hashemi19} and references therein.

In this paper we develop a new algorithm which outperforms the eigenvalue-based methods explored in \cite{Hashemi19} if the degree $n$ is high or the evaluation points $x$ are close to the boundaries of the domain $[-1, 1]$. The new algorithm has two main steps: We first employ the inverse Joukoswki map to convert the problem from the Chebyshev basis into that of Laurent. Then, we apply the interval Horner method to enclose the polynomial in the new basis. Like most of the techniques investigated in \cite{Hashemi19}, the number of basic arithmetic operations involved in the new algorithm is $\mathcal{O}(n)$. 

The most well-known algorithm for evaluation of the polynomial $p$ in {\it floating point arithmetic} is the Clenshaw recurrence \cite{Clenshaw55} which defines the quantities $b_k$ as:
\begin{equation}
\label{clenshawRec:eq}
\left\{ 
\begin{array}{l}
b_{n+2} = b_{n+1} := 0,\\
b_k := 2xb_{k+1} - b_{k+2} + c_k, \ \ k= n, n-1, \dots, 0,
\end{array} 
\right.
\end{equation}
so that $p(x) = b_0 - b_1 x$. On the other hand, an important category of interval arithmetic techniques for enclosing Chebyshev expansions include the parallelepiped and Lohner's QR decomposition methods and two eigenvalue-based algorithms of \cite{Hashemi19}. The basic idea behind all of these techniques is to reformulate the Clenshaw recurrence (\ref{clenshawRec:eq}) in terms of the following discrete dynamical system
\begin{equation}
\label{dynamicalSys:eq}
\hat{b}_{k} = M \hat{b}_{k+1} + \hat{c}_k
\end{equation}
with $M \in \Rz^{2 \times 2}$ and $\hat{b}_{k}, \hat{b}_{k+1}, \hat{c}_{k} \in \Rz^2$ where 
\[
\left. 
\begin{array}{cccccc}
 \underbrace{
 \left( \begin{array}{c}
b_k \\
 b_{k+1}  
 \end{array} \right)}_{\hat{b}_{k}}
 & 
 =
 &
 \underbrace{\left( \begin{array}{cc}
2x & -1  \\
1 & 0  \end{array} \right)}_{M}
  & 
 \underbrace{
\left( \begin{array}{c}
b_{k+1} \\
 b_{k+2} 
 \end{array} \right)}_{\hat{b}_{k+1}}
 & 
 +
 &
\underbrace{
\left( \begin{array}{c}
 c_{k}   \\
0
 \end{array} \right)
 }_{\hat{c}_k}, \ k = n, \dots, 1, 0.
    \end{array} \right.
\]
Unfortunately, the matrix-vector multiplications in \eqnref{dynamicalSys:eq}, when performed in interval arithmetic, cause a severe amount of overestimation called the wrapping effect \cite{Lohner01}.\ To alleviate these overestimations, the eigenvalue-based methods of \cite{Hashemi19} employ the spectral transformation $M = VD V^{-1}$ where
\[
\hspace*{-0.1cm}
\left. \begin{array}{cc}
V  =
\left( \begin{array}{cc}
x +i \sqrt{1-x^2} & x -i \sqrt{1-x^2}  \\
1 & 1  \end{array} \right),
  & 
V^{-1}  =\frac{-i}{2\sqrt{1-x^2}}
\left( \begin{array}{cc}
1 & i \sqrt{1-x^2} - x  \\
-1 &i \sqrt{1-x^2} + x  \end{array} \right),
    \end{array} \right.
\]
\begin{equation}
\label{eigs:eq}
D  = \left( \begin{array}{cc}
x +i \sqrt{1-x^2} & 0  \\
0 & x -i \sqrt{1-x^2}  \end{array} \right).
\end{equation}
Then \eqnref{dynamicalSys:eq} is equivalent to the transformed iteration
\begin{equation}
\label{transformedIter:eq}
\check{b}_{k} = D \check{b}_{k+1} + \check{c}_k,
\end{equation}
where $\check{b}_{k} := V^{-1} \hat{b}_{k}$, and $\check{c}_k := V^{-1} \hat{c}_k$. In practice, interval matrices $\DD \ni D$, $\V \ni V$, and $\I_V \ni V^{-1}$ are used so that rounding errors in the computation of $V$, $D$ and $V^{-1}$ are taken care of. Note that the amount of overestimation, when performing iterations like \eqnref{transformedIter:eq} and \eqnref{dynamicalSys:eq} in interval arithmetic, is governed by the spectral radius of the absolute value of the iteration matrix. As discussed in \cite{Hashemi19}, the reason the transformed iteration \eqnref{transformedIter:eq} outperforms the original iteration \eqnref{dynamicalSys:eq} in interval arithmetic is that $\rho(|D|) = 1$ for every $x \in [-1, 1]$ whereas $1 \leq \rho(|M|) \leq 1 + \sqrt{2}$; see \cite[Fig. 3]{Hashemi19}. On the other hand, it is proved that \cite{Hashemi19}
\[
\kappa_2(V) := \| V^{-1}\| \|V\| =
\left\{ 
\begin{array}{ll} 
\sqrt{\frac{1-x}{1+x}},& -1\leq x \leq 0, \\
\sqrt{\frac{1+x}{1-x}},& 0 \leq x \leq 1.
\end{array} 
\right.
\]
Therefore, the eigenvector matrix $V$ is numerically ill-conditioned for $x \approx \pm 1$. This makes entries of the interval matrix $\I_V$ wide; see e.g., \cite{Revol17} and \cite[p. 346]{Rump10AN}. Consequently, the computed enclosure for $p(x)$ becomes wide. Hence, as discussed in \cite{Hashemi19}, one cannot expect the transformed iteration \eqnref{transformedIter:eq} to give very narrow enclosures at\ignore{points which are close to the boundaries $\pm 1$} $x \approx \pm 1$.\ So, here is the question: can we develop a transformation of \eqnref{dynamicalSys:eq} that takes advantage of the fact that $\rho(|D|) = 1$ while avoiding multiplications by the eigenvector matrix $V$ and its inverse? It turns out that the following method is what we are looking for.

\section{The Laurent-Horner method}The Joukowski map $x := J(z) = \frac{1}{2}(z + \frac{1}{z})$
\ignore{\[
x := J(z) = \frac{1}{2}(z + \frac{1}{z})
\]}
is a popular conformal map in approximation theory and complex analysis. It transforms origin-centered circles to ellipses with foci at $\{-1, 1\}$ which are known as {\em Bernstein ellipses}. In particular, it maps the unit circle to the unit interval $[-1, 1]$. Since $J(z) = J(z^{-1})$, there is a 2-to-1 correspondence between $z$ on the unit circle and $x$ on the real interval $[-1, 1]$. In other words, the quadratic equation corresponding to the map has two solutions $z = x \pm i \sqrt{1-x^2}$. Notice that these are the eigenvalues of $M$ in
\eqnref{dynamicalSys:eq} as can be observed from the spectral decomposition \eqnref{eigs:eq}. 

Another basic fact that is important for the development of the new method is that the $k$-th Chebyshev polynomial is the real part of the function $z^k$ on the unit circle \cite{ATAP}, i.e.,
\[
T_k(x) = \frac{z^k + z^{-k}}{2} = \frac{z^k + \bar{z}^k}{2},
\]
where $\bar{z}$ denotes the complex conjugate. Hence, the Chebyshev expansion \eqnref{ChebExpand: eq} can be converted to the following {\em Laurent} polynomial
\begin{equation}
\label{LaurExpand: eq}
p(x) = \frac{1}{2}\sum_{k=0}^n c_k (z^k + z^{-k}) = \mbox{ real }\Big(\sum_{k=0}^n c_k z^k\Big).
\end{equation}
To take care of rounding errors in the conversion from $x$ to $z$, we compute an interval $\z$ containing the exact value of $z$. 

The second step of our enclosure method simply applies the Horner's rule in interval arithmetic to \eqnref{LaurExpand: eq}; see e.g., \cite{Ceberio02}.\ The interval Horner method is a straightforward extension of the standard Horner's nested multiplication form to interval arithmetic and can be used to enclose the range of polynomials with a linear complexity in terms of its degree. 

\ignore{Do we want to say more about the interval Horner method especially its optimality in certain cases?
BH: No need to mention optimality conditions of Stahl as they only apply to evaluation at a REAL interval. However, it might be good to explain \cite{Ceberio02}.
  Can we relate them more specifically to our problem here? Horner's rule is a special case of the Clenshaw recurrence evaluation of monomials in the monomial basis. }

While {\em both} the existing eigenvalue-based methods and the Laurent-Horner method inherit possible ill-conditioning of the transformation from $x$ to $z$, computations in the Laurent-Horner method do not involve the matrices $V$ and $V^{-1}$. Therefore, the new method might be considered as an eigenvector-avoiding ``spectral'' transformation of \eqnref{dynamicalSys:eq} and can be expected to outperform the eigenvalue-based methods of \cite{Hashemi19} especially at $x \approx \pm 1$. 

\section{Numerical experiments}
To compare the new algorithm with the older ones, we illustrate the time needed (in seconds) together with $\rad \p(x)$; the radius of the computed enclosures. The radii are employed also for obtaining the average number of correct digits for enclosures over all points $x$ computed as ${\rm mean} ( -\log_{10} \big(\rad \p(x) \big) )$.
Our numerical results are generated using INTLAB \cite{Rump99b}.

\begin{figure}[!h]
\begin{center}
\includegraphics[width=0.8\textwidth, height=7cm]{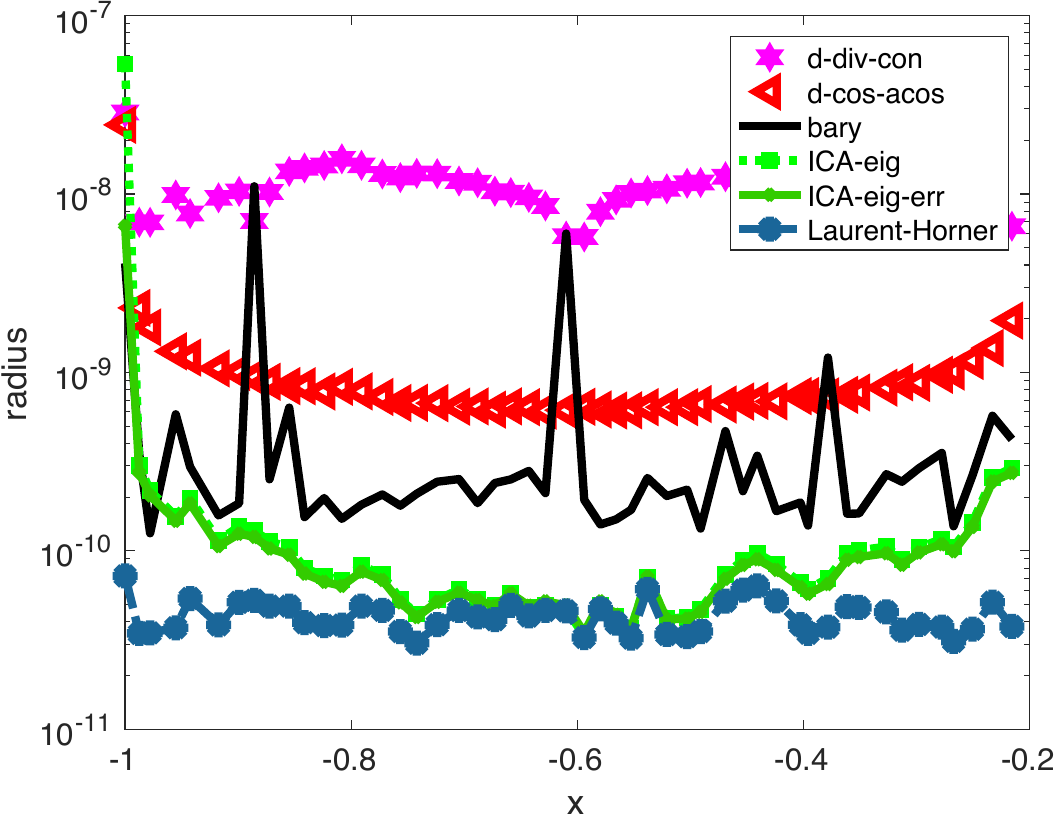}
\caption{Radius of enclosures for a degree 9150 interval polynomial at 1000 intervals. Only 50 points are depicted to make the curves easier to distinguish. \label{ex9_1:fig}}
\end{center}
\end{figure}


\begin{figure}
\begin{center}
\includegraphics[width=0.8\textwidth, height=6.8cm]{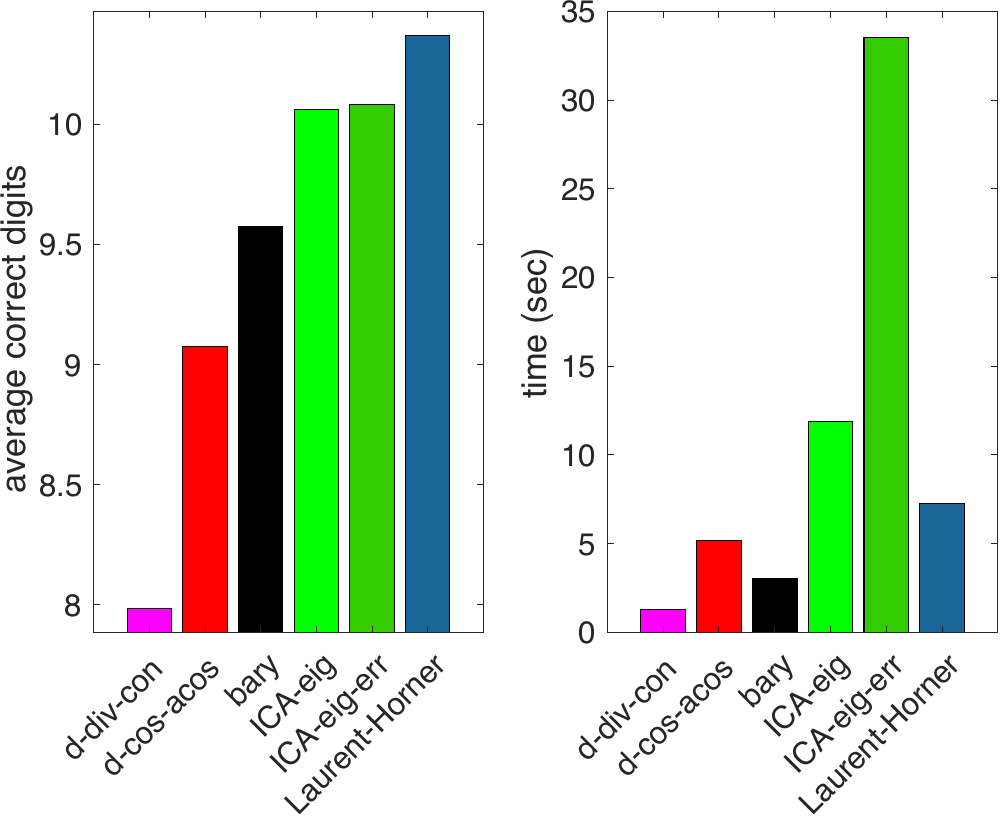}
\caption{Average number of correct digits (left) and computing time (right) of different methods for bounding the range of a Chebyshev expansion of degree 9150 at 1000 points. \label{ex9_2:fig}}
\end{center}
\end{figure}

We consider a Chebyshev expansion of degree $9150$ corresponding to the random smooth function \texttt{randfun(0.0007)} \cite{Filip18} where its real coefficients are inflated to be intervals of a radius of about $2\times 10^{-15}$. The same Chebyshev expansion is used in Example 6.6 of \cite{Hashemi19}. We compute enclosures for the value of the interval polynomial at $1000$ random intervals $\x$ whose radii are again of the order of $10^{-15}$. Figures \ref{ex9_1:fig} and \ref{ex9_2:fig} contain our results for six methods. Here, \texttt{ICA-eig} and \texttt{ICA-eig-err} denote the two eigenvalue-based methods of \cite{Hashemi19} which rely on \eqnref{eigs:eq} and \eqnref{transformedIter:eq}. Also, \texttt{d-cos-acos} and \texttt{d-div-con} denote two direct methods which typically give narrow enclosures for low-degree Chebyshev expansions. Moreover, \texttt{bary} refers to the extension of barycentric formula to interval arithmetic.



It was shown in \cite{Hashemi19} that in the case of high-degree Chebyshev expansions, \texttt{ICA-eig} and \texttt{ICA-eig-err} typically give narrowest enclosures among the vectorized techniques. 
Nevertheless, we observe that in this example the Laurent-Horner method not only outperforms those techniques with respect to speed, but also computes narrowest enclosures among {\it all} the techniques. 

\ignore{We end this paper with a note concerning the conversion of real evaluation points $x$ into complex intervals $\z$ in the first step of the Laurent-Horner method. This might seem a bad idea both with respect to the slower speed and more overestimation in complex interval arithmetic as opposed to its real counterpart. However, this is not a concern in INTLAB where working with complex intervals does not cause much penalty. The reason is that arithmetic operations on both real and complex intervals as implemented in INTLAB use the same midpoint-radius representation; see \cite{Rump99a} for details.}

In a second variant of the above experiment, narrowest-possible enclosures for $c_k$ and $x$ are used instead of those with a width of order $10^{-15}$. This time, the most accurate methods {\em on average} are \texttt{ICA-eig-err} and Laurent-Horner, while \texttt{ICA-eig-err} is three times slower than the Laurent-Horner.

Let us end this paper with a note concerning the conversion from real to complex interval arithmetic in the new method as well as in \texttt{ICA-eig} and \texttt{ICA-eig-err}. The conversion has a speed penalty observed for the three slower methods in Figure \ref{ex9_2:fig} (right). However, moving to the complex plane, while avoiding wrappings caused by the eignevector transformations 
are the main reasons the new method gives the narrowest enclosures. Avoiding these multiplications also makes the new method the fastest among those that use complex arithmetic.\ Note also that INTLAB employs midpoint-radius representation in its implementation of complex as well as real machine interval arithmetic \cite{Rump99a}.



\paragraph{Acknowledgement} We thank an anonymous referee for valuable comments. Jared L. Aurentz has received financial support through the
Postdoctoral Junior Leader Fellowship Programme from ``la Caixa'' Banking Foundation. The second author was supported by the Iranian National Science Foundation (INSF) under grant No. 98012590.

\bibliographystyle{elsarticle-num}
\bibliography{paperbib}

\end{document}